\documentclass[xypic,amssymb,rawfonts,amsthm,amsmath,amsfonts,syntonly,11pt]{amsart}

\usepackage{amssymb}
\usepackage{hyperref}
\usepackage{mathrsfs}

\newcounter{thecounter}
\numberwithin{thecounter}{section}

\newtheorem{thm}[thecounter]{Theorem}

\theoremstyle{definition}

\numberwithin{equation}{section}

\newcommand{\map}{\operatorname{map}}

\newcommand{\aut}{\operatorname{aut}}

\newcommand{\Aut}{\operatorname{Aut}}

\newcommand{\Out}{\operatorname{Out}}
\newcommand{\End}{\operatorname{End}}

\newcommand{\im}{\operatorname{im}}

\newcommand{\hocolim}{\operatorname{hocolim}}

\newcommand{\cC}{{\mathcal{C}}}

\newcommand{\bO}{{\mathbf O}}

\newcommand{\pcom}{\hat{{}_p}}

\newcommand{\twocom}{\hat{{}_2}}

\newcommand{\PU}{\operatorname{PU}}

\newcommand{\SU}{\operatorname{SU}}
\newcommand{\Sp}{\operatorname{Sp}}

\newcommand{\Spin}{\operatorname{Spin}}

\newcommand{\F}{{\mathbb {F}}}
\newcommand{\Z}{{\mathbb {Z}}}

\newcommand{\D}{{\mathbf D}}

\newcommand{\A}{{\mathbf A}}

\newcommand{\W}{{\mathcal{W}}}
\newcommand{\N}{{\mathcal{N}}}

\newcommand{\cZ}{{\mathcal{Z}}}

\newcommand{\beq}{\begin{eqnarray*}}
\newcommand{\eeq}{\end{eqnarray*}}

\newcommand{\tuborg}{\left\{\begin{array}{ll}}
\newcommand{\sluttuborg}{\end{array}\right.}

\newfont{\bm}{msbm10}

\input xypic
\xyoption{all}
\CompileMatrices       

\renewcommand{\phi}{\varphi}
\begin{document}
\title{The classification of $2$-compact groups} 


\thanks{Partially supported by NSF grant DMS-0354633}

\author[J. Grodal (joint with K. Andersen)]{Jesper Grodal \\ (joint work with Kasper K.\ S.\ Andersen)}



\address{Department of Mathematics, University of Aarhus, DK-8000 {{\AA}rhus C}, Denmark}

\email{kksa@math.ku.dk}

\address{Department of Mathematics, University of Chicago, Chicago, IL 60637, USA}
\email{jg@math.uchicago.edu}
\begin{abstract} 
In this talk I'll announce and explain a
proof of the classification of $2$-compact groups, joint with K.\ Andersen, hence completing the classification of
p-compact groups at all primes $p$. A $p$-compact group, as introduced by Dwyer-Wilkerson, is
a homotopy theoretic version of a compact Lie group, but with all its structure
concentrated at a single prime $p$. Our classification states that there is a
$1$-$1$-correspondence between connected $2$-compact groups and root data over the $2$-adic
integers (which will be defined in the talk). As a consequence we get the conjecture that
every connected $2$-compact group is isomorphic to a product of the $2$-completion of a
compact Lie group and copies of the exotic $2$-compact group $DI(4)$, constructed by
Dwyer-Wilkerson. The major new input in the proof over the proof at odd primes (due to
Andersen-Grodal-M{\o}ller-Viruel) is a thorough analysis of the concept of a root datum for
$2$-compact groups and its relationship with the maximal torus normalizer. With these tools
in place we are able to produce a proof which to a large extent avoids case-by-case
considerations.
\end{abstract}

\maketitle



This is a summary of my lecture given at the ``Conference on Pure and Applied Topology'', Isle of Skye, Friday June 24, 2005. I announced and sketched a proof of the classification of $2$-compact groups, in joint work with Kasper Andersen. This will appear in the papers \cite{AG05automorphisms} and \cite{AG05classification}.

Recall that a $p$-compact group is a triple $(X,BX,e: X \xrightarrow{\simeq} \Omega BX)$ where $BX$ is a pointed, connected, $p$-complete space of the homotopy type of a $CW$-complex, $X$ satisfies that $H^*(X;\F_p)$ is finite over $\F_p$, and $e$ is a homotopy equivalence.  They are homotopy theoretic analogs of compact Lie groups, and were introduced by Dwyer-Wilkerson in \cite{DW94}.
 A $p$-compact group is said to be connected if $X$ is a connected space.

Our main theorem is the following
\begin{thm}[\cite{AG05classification}] Let $(X,BX, e: X \xrightarrow{\simeq} \Omega BX)$ be a connected $2$-compact group.
Then
$$BX \simeq BG\twocom \times BDI(4)^s$$
where $BG\twocom$ is the $2$-completion of a connected compact Lie group $G$, and $BDI(4)$ is the classifying space of the exotic $2$-compact group $DI(4)$ constructed by Dwyer-Wilkerson in \cite{DW93}, $s \geq 0$.
\end{thm}
A corresponding statement for odd primes was proved by the authors together with M{\o}ller and Viruel in \cite{AGMV03}. Partial results for $p=2$ have been obtained by Dwyer-Miller-Wilkerson, Notbohm, Viruel, M{\o}ller, and Vavpeti\v c-Viruel. Our proof is a self-contained induction.

There is a better, more precise, formulation of our theorem which both makes it clear why it is the correct $2$-local version of the classification of compact Lie groups and suggests a possible strategy of proof. It is based on the notion of a root datum over the $p$-adic integers $\Z_p$, which we now introduce.

For $R$ a principal ideal domain, an $R$-root datum $\D$ is defined to
be a triple $(W,L,\{Rb_\sigma\})$, where $L$ is a free
$R$-module of finite rank, $W \subseteq \Aut_R(L)$ is a finite
subgroup generated by reflections (i.e., elements $\sigma$ such that
$1-\sigma \in \End_{R}(L)$ has rank one), and $\{Rb_\sigma\}$
is a collection of rank one submodules of $L$, indexed by the
reflections $\sigma$ in $W$, satisfying
$$
\im(1-\sigma) \subseteq Rb_\sigma \subseteq \ker(
\sum_{i=0}^{|\sigma|-1} \sigma^i) \mbox{ and } w (Rb_\sigma) =
Rb_{w\sigma w^{-1}} \mbox{ for all } w \in W
$$
The element $b_\sigma \in L$, called the {\em coroot} corresponding to
$\sigma$, is determined up to a unit in $R$. Together with $\sigma$ it determines a root $\beta_\sigma: L \to R$ via the formula
$$
\sigma(x) = x + \beta_\sigma(x)b_\sigma
$$
If $R = \Z$ then there is a $1$-$1$-correspondence between $\Z$-root
data and classically defined root data, by to $(W,L,\{\Z
b_\sigma\})$ associating 
 $(L,L^*,\{\pm b_\sigma \},\{\pm \beta_\sigma \})$; 
see \cite[Prop.~2.16]{DW05}.
If $R
= \Z$ or $\Z_p$, instead of the collection $\{Rb_\sigma\}$ one can
equivalently consider their span, the {\em coroot lattice}, $L_0 =
+_\sigma Rb_\sigma \subseteq L$, and this was the definition given in
\cite[\S~1]{AGMV03}, under the name ``$R$-reflection datum''.
Two $R$-root data $\D = (W,L,\{Rb_\sigma\})$ and $\D' =
(W',L',\{Rb'_\sigma\})$ are said to be isomorphic if there exists an
isomorphism $\phi: L \to L'$ such that $\phi W\phi^{-1} = W'$ and
$\phi(Rb_\sigma) = Rb'_{\phi \sigma \phi^{-1}}$. In particular the
automorphism group is given by $\Aut(\D) = \{ \phi \in
N_{\Aut_R(L)}(W) | \phi(Rb_\sigma) = Rb_{\phi(\sigma)} \}$ and we
define the outer automorphism group as $\Out(\D) = \Aut(\D)/W$.

We now explain how to associate a $\Z_p$-root datum to a $p$-compact group. By a theorem of Dwyer-Wilkerson any $p$-compact group $(X,BX,e)$ has a maximal torus, which is a map $i: BT = (BS^1\pcom)^r \to BX$ satisfying that the fiber has finite $\F_p$-cohomology and non-trivial Euler characteristic. Replacing $i$ by an equivalent fibration, we define the 
 {\em Weyl space} ${\W}_X(T)$ as the topological monoid of
self-maps $BT \to BT$ over $i$. The {\em Weyl group} is defined as
$W_X(T) = \pi_0({\W}_X(T))$ and the classifying space of the maximal torus normalizer is defined as $B\N_X(T) = BT_{h\W_X(T)}$. Now, by definition, $W_X$ acts on $L = \pi_2(BT)$ and it is a theorem of Dwyer-Wilkerson that, if $X$ is connected, this gives a faithful representation of $W_X$ in $\Aut_{\Z_p}(L)$ as a finite $\Z_p$-reflection group. There is also an easy formula for the $b_\sigma$ in terms of the maximal torus normalizer $\N_X$, for which we refer to \cite{DW05}, \cite{AG05automorphisms}, or \cite{AG05classification}. This explains how to associate a $\Z_p$-root datum $\D_X$ to a connected $p$-compact group $X$. (In the case where $p$ is odd the root datum is in fact completely determined by the finite $\Z_p$-reflection group $(W,L)$, which explains the formulation of our classification theorem with M{\o}ller and Viruel in \cite{AGMV03}.)

We are now ready to state the precise version of our main theorem.
\begin{thm}[\cite{AG05classification}] \label{classthmroot} The assignment which to a connected $2$-compact group $X$ associates its $\Z_2$-root datum $\D_X$ root gives a one-to-one correspondence between connected $2$-compact groups and $\Z_2$-root data.  Furthermore the map $\Phi: \pi_0(\Aut(BX))  \to \Out(\D_X)$ is an isomorphism, and $B\Aut(BX)$ is the unique total space of a split fibration
$$ B^2\cZ(\D) \to B\Aut(BX) \to B\Out(\D_X)$$
\end{thm}
Here $\cZ(\D)$ is the center of the root datum $\D$, defined so as to agree with the formula for the center of a $p$-compact group given in \cite{dw:center}, and $B\Aut(BX)$ denotes the classifying space of the topological monoid of self-homotopy equivalences of $BX$.

We remark that since we have completely described connected $p$-compact groups {\em as well as their space of automorphisms}, we also get a classification for non-connected $2$-compact groups. (The situation here is totally analogous to the case of compact Lie groups.) If $X$ is a connected $p$-compact group with a given root datum $\D$ and $\pi$ is a given $p$-group, then the $p$-compact groups with connected component isomorphic to $X$ and component group $\pi$ are in $1$-$1$-correspondence with the $\Out(\pi)$-orbits of the set $[B\pi,B\Aut(BX)]$. Since both $X$ and and $B\Aut(BX)$ is completely described in terms of $\D_X$ by the main theorem, this gives a classification also in the non-connected case. We refer to \cite{AG05classification} for a more detailed description.
\smallskip

The main theorem has a number of corollaries. The most important is perhaps that it gives a proof of the maximal torus conjecture, giving a purely homotopy theoretic characterization of compact Lie groups amongst finite loop spaces.

\begin{thm}[Maximal torus conjecture \cite{AG05classification}]  The classifying space functor, which to a compact Lie group $G$ associates the finite loop space $(G,BG,e: G \xrightarrow{\simeq} \Omega BG)$ gives a one-to-one correspondence between compact Lie groups and finite loop spaces with a maximal torus. 
(Moreover, if $G$ is connected we have a split fibration
$B^2\cZ(\D_G) \to B\Aut(BG) \to  B\Out(\D_G)$.)
\end{thm}
The fact that the functor ``$B$'' is faithful was already known by work of Notbohm, M{\o}ller, and Osse and the statement about the space $B\Aut(BG)$ follows easily from earlier work of Jackowski-McClure-Oliver and Dwyer-Wilkerson. The new, and a priori quite surprising, result here is the statement that if a finite loop space has a maximal torus, then it has to come from a compact Lie group.

Another application of the classification of $2$-compact groups is to give an answer to the so-called Steenrod problem for $p=2$ (see \cite{steenrod61} and \cite{steenrod71}), which asks which graded polynomial algebras can occur as the mod $2$ cohomology ring of a space? Steenrod's problem was solved for $p$ ``large enough'' by Adams-Wilkerson \cite{AW80} and for all odd primes by Notbohm \cite{notbohm99} using a partial classification of $p$-compact groups, $p$ odd.

\begin{thm}[Steenrod's problem for $p=2$ \cite{AG05classification}] Suppose that $P^*$ is a graded polynomial algebra over $\F_2$ in finitely many variables. If $P^*$ occurs as $H^*(Y;\F_2)$ for some space $Y$ (no assumptions), then $P^*$ is isomorphic, as a graded algebra, to $$H^*(BG;\F_2) \otimes H^*(BDI(4);\F_2)^{\otimes s} \otimes Q^*$$ where $G$ is a connected semi-simple Lie group and $Q^*$ is a polynomial ring with generators in degrees one and two.

In particular if $P^*$ is assumed to have generators in degree $\geq 3$  then $G$ has to be simply connected and $P^*$ is a tensor product of the following graded algebras:

\begin{center}
\begin{tabular}{lr}
$\F_2[x_4, x_6, \ldots, x_{2n}]$ & $(\SU(n))$\\
$\F_2[x_4,x_8, \mbox{ }\ldots, x_{4n}]$ & $(\Sp(n))$\\
$\F_2[x_4, x_6, x_7, x_8]$ & $(\Spin(7))$\\
$\F_2[x_4, x_6, x_7, x_8, x'_8]$ & $(\Spin(8))$\\
$\F_2[x_4, x_6, x_7, x_8, x_{16}]$ & $(\Spin(9))$\\

$\F_2[x_4,x_{6},x_{7}]$ & $(G_2)$\\
$\F_2[x_4, x_6, x_7, x_{16}, x_{24}]$ & $(F_4)$\\
$\F_2[x_8,x_{12},x_{14},x_{15}]$ & $(DI(4))$\\
\end{tabular}
\end{center}
\end{thm}
It seems reasonable that one can in fact list all polynomial rings which occur as $H^*(BG;\F_2)$ for $G$ semi-simple, although we have not been able to locate such a list in the literature; for $G$ simple a list can be found in \cite{kono75}.

We finally point out that many classical theorems from Lie theory by Borel, Bott, Demazure, and others also carry over to $2$-compact groups via the classification. Applications of this type were already pointed out in \cite{AGMV03}, to which we refer.

\smallskip

We end this short summary by saying a few words about the proofs. Our strategy is an elaboration of the strategy we pursued with M{\o}ller and Viruel in \cite{AGMV03}, but with several new additions, most importantly a development and utilization of the theory of $p$-adic root data. The sketch below is intended to sum up the main steps without introducing too much notation. We focus on the uniqueness statement, that two connected $p$-compact groups with isomorphic root data are isomorphic.

It is a recent theorem of Dwyer-Wilkerson \cite{DW05}, extending old work of Tits, that the maximal torus normalizer $\N_X$ in a $2$-compact group $X$ can be recovered from the root datum $\D_X$. However, $\N_X$ has a larger group of automorphisms than $X$, which means that it is not the right classifying invariant. 

This problem can be overcome, by also keeping track of certain ``root subgroups'', a machinery which we set up in \cite{AG05automorphisms}.
In a given maximal torus normalizer $\N$, for each reflection $\sigma$ and coroot
$b_\sigma$
one can algebraically construct a root subgroup $\N_\sigma$ of
$\N$. 
In the
setting of algebraic groups, this ``root subgroup'' $\N_\sigma$ will
be the maximal torus normalizer of $\langle U_\alpha,
U_{-\alpha}\rangle$, where $U_\alpha$ is the root subgroup in the
sense of algebraic groups corresponding to the root $\alpha$ dual to the coroot
$b_\sigma$. 

\begin{thm}[\cite{AG05automorphisms}] \label{autothm} Let $X$ be a connected $2$-compact group.
The canonical map $\Out(BX) \to \Out(\D_X)$ factors 
$$\Phi: \Out(BX) \to \Out(B\N,\{B\N_\sigma\}) \xrightarrow{\simeq} \Out(\D_X)$$
and we likewise have a canonical map
$$\Phi: B\Aut(BX) \to B\aut(\D_X)$$ 
where $B\aut(\D_X)$ is a certain space defined in \cite{AG05automorphisms} refining $B\Aut(B\N,\{B\N_\sigma\})$.
\end{thm}
These maps will be shown to be isomorphisms inductively, as part of proving Theorem~\ref{classthmroot}.
Since the pair $(\N,\{\N_\sigma\})$ carries the same structure and automorphisms as $\D$, but is closer to $X$, it is more useful than $\D$ for classification purposes.

Assume that we have two connected $p$-compact groups $X$ and $X'$
with the same root datum $\D$ and hence the same maximal torus
normalizer and root subgroups $(\N,\{\N_\sigma\})$, and assume by
induction that Theorem~\ref{classthmroot} is true for all
connected $p$-compact groups of smaller cohomological dimension. The
first step is to observe that one can reduce to the case where $X$ and
$X'$ are simple and center-free, mimicking the proofs in \cite{AGMV03} but everywhere keeping track of the root subgroups. So we can assume that we are in the following situation
$$ \xymatrix{ & (B\N,\{B\N_\sigma\}_\sigma) \ar[dr] \ar[dl] &\\
                BX \ar@{-->}[rr]&& BX'}$$
where the dotted arrow is the one that we want to construct.

Using that every element of order $p$ can be conjugated into the maximal
torus, uniquely up to conjugation in $\N$, we can by induction construct the dotted arrow on the connected component of the centralizer of every element $\nu: B\Z/p \to BX$ of order $p$ in $X$. 
$$ \xymatrix{ & (B\cC_\N(\nu)_1,\{B{\cC_\N(\nu)_1}_\sigma\}_\sigma) \ar[dr] \ar[dl] &\\
                B\cC_X(\nu)_1 \ar[rr]& & B\cC_{X'}(\nu)_1 }$$
Using that we have control of the whole
{\em space} of self-equivalences of $B\cC_X(\nu)_1$ by
Theorem~\ref{autothm} and induction, one sees that this map is
equivariant with respect to the component group and hence extends
uniquely to the whole centralizer $B\cC_X(\nu)$.

Now for a general elementary abelian subgroup $\nu: E \to X$ of $X$ we can pick an element of order $p$ in $E$, and we get via restriction a map
$$B\cC_X(\nu) \to B\cC_{X}(\Z/p) \to B\cC_{X'}(\Z/p) \to BX'$$
To make sure that these maps are chosen in a compatible way, one has
to show that this map does not depend on the choice of rank one
subgroup of $E$. We developed some techniques in \cite{AGMV03} for
handling this kind of situation. Using easy (one page)
 case-by-case arguments on the level of
$2$-adic root data, we are able to verify that the assumptions
of \cite{AGMV03} are always satisfied. 

Since these maps from centralizers are chosen in a compatible way, they combine to form an element
$$\vartheta \in \lim_{\nu \in \A(X)}\hspace{-10pt}{}^0[B\cC_X(\nu),BX']$$
where $\A(X)$ is the Quillen category of $X$. Furthermore the construction of this element allows one to show that $X$ and $X'$ have the same $p$-fusion and in particular the same maximal torus $p$-normalizer and the same $\Z_p$-cohomology. We are left with a rigidification question.

One approach could be to show that the relevant obstruction groups
$$\lim_{\vartheta \in \A(X)}\hspace{-10pt}{}^*\pi_*(\map(B\cC_X(\nu),BX')_{[\vartheta]})$$
vanish, which would show that our diagram rigidifies and produces an equivalence
$$BX \simeq \hocolim_{\nu \in \A(X)}B\cC_X(\nu) \xrightarrow{\simeq} BX'$$
It is immediate to see that this is the case for many groups $X$ including $DI(4)$, and it was the approach taken for $p$ odd in \cite{AGMV03}. There we for instance also show that it holds for $X = \PU(n)\twocom$, by the same calculation as for $p$ odd, but it is more difficult to verify e.g., for the exceptional groups.

We however take a different approach: Since $DI(4)$ is easily dealt with, we can assume that $BX = BG\twocom$ for some simple, center-free Lie group. For each $p$-radical subgroup $P$ of $G$ we have a map $BP \to BC_G({}_pZ(P))\pcom \to BX'$. Since $X$ and $X'$ have the same $\Z_p$-homology, and in particular the same fundamental group, we can lift this map to a map $B\tilde P \to B\tilde X'$, where $\tilde P$ is the preimage of $P$ in $\tilde G$, the universal cover of $G$, and $\tilde X'$ is the universal cover of $X'$.

One now verifies that these maps assemble to give an element in 
$$ \lim_{\tilde G/\tilde P \in \bO_p^r(\tilde G)}\hspace{-15pt}{}^0[B\tilde P,B\tilde X']$$
where $\bO_p^r(\tilde G)$ is the full subcategory of the $p$-orbit category of $\tilde G$ with objects $\tilde G/\tilde P$ for $\tilde P$ a $p$-radical subgroup of $\tilde G$.
The obstructions to rigidifying this to get a map 

$$ B\tilde G \simeq \hocolim_{\tilde G/\tilde P \in \bO_p^r(\tilde G)} B\tilde P \to B\tilde X'$$
can be shown to lie in groups which identify with
$$\lim_{\tilde G/\tilde P \in \bO_p^r(\tilde G)}\hspace{-15pt}{}^*\pi_*(Z(\tilde P))$$
These groups have been shown to identically vanish by earlier
work of Jackowski-McClure-Oliver \cite{JMO92}. Hence the map exists and it now follows easily by the construction that it is an equivalence. Passing to a quotient we get the wanted equivalence $BG \to
BX'$, finishing the proof of the theorem.


\bibliographystyle{plain}
\bibliography{../pcg/poddclassification}
\end{document}